\documentclass[12pt]{amsart}
\usepackage{latexsym,amsmath, amscd, amssymb, amsthm, euscript}
\newtheorem{thm}[subsection]{Theorem}

\newtheorem{lem}[subsection]{Lemma}

\newtheorem{prop}[subsection]{Proposition}
\theoremstyle{definition}

\newtheorem{defn}[subsection]{Definition}

\theoremstyle{remark}

\theoremstyle{remark}

\newtheorem{rem}[subsection]{Remark}

\numberwithin{equation}{subsection}

\newcommand{\Q}{\mathbb{Q}}

\newcommand{\mc}{\mathcal }

\newcommand{\scr}[1]{\mathbf{\EuScript{#1}}}

\newcommand{\Z}{\mathbb{Z}}

\newcommand{\Sp}{\text{\rm Spec}}

\newcommand{\ru}[1]{\ulcorner #1 \urcorner}
\newcommand{\fp}[1]{\langle #1 \rangle}
\begin{document}

\title{Kawamata--Viehweg vanishing as Kodaira vanishing for stacks}%
\author{Kenji Matsuki and Martin Olsson}%

\begin{abstract} We associate to a pair $(X,D)$, consisting of a smooth scheme
with a divisor $D\in \text{Div}(X)\otimes \mathbb{Q}$ whose
support is a divisor with normal crossings, a canonical
Deligne--Mumford stack over $X$ on which $D$ becomes integral. We
then reinterpret the Kawamata--Viehweg vanishing theorem as
Kodaira vanishing for stacks.
\end{abstract}

\maketitle

\date{\today}%

\section{Introduction}
Let $k$ be a field, and suppose given a pair $(X, D)$, where $X/k$
is a smooth variety and $D\in \text{Div}(X)\otimes \Q$ is a
divisor whose support has normal crossings. Suppose further that
there exists an integer $n$, prime to $\text{char}(k)$, for which
$nD$ is integral. Our aim in this note is to explain how one can
associate to $(X, D)$ a Deligne--Mumford stack $\scr X/X$, which
in a precise sense is the ``minimal covering'' of $X$ on which $D$
becomes integral (see (\ref{thm3}) for the precise properties of
$\scr X$). In addition, we explain, using the stack $\scr X$, how
the Kawamata--Viehweg vanishing theorem (\cite{E-V}, \cite{Kaw},
\cite{Vie}) can naturally be interpreted as an application of
Kodaira vanishing for stacks.

The note is organized as follows.  In section 2 we prove a Kodaira
vanishing theorem for Deligne--Mumford stacks.  A Hodge theoretic
approach to the  characteristic $0$ version of this theorem had
previously been considered by Starr (private communication).  In
section 3 we make a few basic observations concerning log
structures in the sense of Fontaine and Illusie which will be
needed in the construction of $\scr X$. In section 4 we construct
the stack $\scr X$ and describe its basic properties. Then in
section 5, we discuss the Kawamata--Viehweg vanishing theorem.

\noindent \emph{Acknowledgements:} We would like to thank D.
Arapura, Y. Kawamata,  T. Pantev, and J. Starr
 for their comments
and help.

\section{Kodaira vanishing for stacks}

\begin{thm}\label{thm2} Let $k$ be a field, and let $\scr X/k$ be a smooth proper tame  Deligne--Mumford
stack  of dimension $d$ with projective  coarse moduli space $\pi
:\scr X\rightarrow X$.  Suppose that $\pi $ is flat  and let $L$
be an invertible sheaf on $\scr X$ such that some power of $L$
descends to an ample invertible sheaf on $X$.

\noindent {\rm (i)}. If $\text{\rm char}(k) = 0$, then $H^j(\scr
X, \Omega ^i_{\scr X/k}\otimes L^{-1}) = 0$ for $i+j<d$.

\noindent {\rm (ii)}.  If $\text{\rm char}(k) = p> 0$, $k$ is
perfect, and $\scr X$ admits a smooth lifting to $W_2(k)$, then
$H^j(\scr X, \Omega ^i_{\scr X/k}\otimes L^{-1}) = 0$ for
$i+j<\text{\rm inf} (d, p).$
\end{thm}

Theorem (\ref{thm2}) is proven by the same argument used in the
proof by Deligne, Illusie, and Raynaud of the Kodaira vanishing
theorem (\cite{D-I}, 2.8). We indicate here only the argument
needed for the reduction from characteristic $0$ to characteristic
$p$, and a basic proposition. The rest of the proof can then be
transcribed word for word from (loc. cit.).

\noindent \emph{Reduction to positive characteristic.}

This is fairly standard using the classical argument for schemes
(\cite{EGA}. IV.8.10.5 (xii)).  However, we explain the reduction since it uses a
non--trivial result about stacks:  Chow's lemma for
Deligne--Mumford stacks (\cite{La}, 16.6.1). So consider $\scr
X/k$ where $k$ is a field of characteristic $0$, and choose a
presentation for $\scr X$. That is, suppose $\scr X$ is defined by
affine schemes $U$ and $R$ together with maps $s, b:R\rightarrow
U$ (source and target), $\epsilon :U\rightarrow R$ (identity),
$i:R\rightarrow R$ (inverse), $m:R\times _{s, U, b}R\rightarrow R
$ (composition), such that the finitely many conditions for these
to define an \'etale groupoid in schemes are satisfied (\cite{La},
2.4.3). Since all these schemes are affine of finite type over
$k$, we can find a finitely generated $\mathbb{Z}$--algebra
$A\subset k$, and an \'etale groupoid in schemes $(\widetilde U,
\widetilde R, \tilde s, \tilde b, , \tilde i, \tilde m)$ over $A$
inducing our groupoid over $k$.  Thus we obtain a Deligne--Mumford
stack $\widetilde {\scr X}/A$ inducing $\scr X$.  Now by
(\cite{La}, 16.6.1), there exists a projective scheme $Y/k$
together with a proper surjective morphism $f:Y\rightarrow \scr
X$.  After localizing on $\Sp (A)$ (and using (\cite{La}, 4.18)),
we can therefore find a a projective scheme $\widetilde Y$ and a
map $\tilde f:\widetilde Y\rightarrow \widetilde {\scr X}$.  Now
the condition that the map $\tilde f$ is proper and surjective is
equivalent to the condition that the induced map of schemes
$\widetilde Y\times _{\widetilde {\scr X}}\widetilde U\rightarrow
\widetilde U$ is proper and surjective, and hence after localizing
on $\Sp (A)$ we can assume that $\tilde f$ is proper and
surjective.  From this it follows that $\widetilde {\scr X}$ is
proper. Let $\tilde \pi :\widetilde {\scr X}\rightarrow \widetilde
X$ be its coarse moduli space. There is a maximal open substack of
$\widetilde {\scr X}$ which is smooth over $\Sp (A)$ and flat over
$\widetilde X$, and since $\widetilde {\scr X}$ is proper over $A$
this implies that after localizing on $\Sp (A),$ we can assume
that $\widetilde {\scr X}/A$ is tame and smooth and that $\tilde
\pi $ is flat.  Now giving the invertible sheaf $L$ on $\scr X$ is
equivalent to giving a morphism $\scr X\rightarrow B\mathbb{G}_m$.
Hence by (\cite{La}, 4.18), we can after further localization on
$\Sp (A)$ find a sheaf $\widetilde L$ on $\widetilde {\scr X}$
inducing $L$ such that some tensor power of $\widetilde L$
descends to an ample invertible sheaf on $\widetilde X$.

\begin{lem}\label{2:lem} Let $\scr X/S$ be a tame Deligne--Mumford stack of
finite type over a noetherian scheme $S$, and suppose the map $\pi
:\scr X\rightarrow X$ from $\scr X$ to its coarse moduli space is
flat.  Then

\noindent
 {\rm (i).}  If $F$ is a coherent sheaf on $\scr X$,
then $R^i\pi _*F = 0$ for $i>0$, and the formation of $\pi _*F$ is
compatible with arbitrary base change $S'\rightarrow S$.

\noindent {\rm (ii).} If $E$ is a locally free sheaf on $\scr X$,
then $\pi _*E$ is locally free on $X$.

\noindent {\rm (iii).} $X$ is Cohen--Macaulay.
\end{lem}
\begin{proof}
It is well known that \'etale locally on $X$, the stack $\scr X$
is isomorphic to  $[U/\Gamma ]$, where $U=\Sp (R)$ is an affine
scheme and $\Gamma $ is a finite group of order invertible in $k$
acting on $U$ (see for example (\cite{A-V}, 2.2.3)). Moreover, in
this situation $X = \Sp (R^\Gamma )$. Now if $M$ is any $R^\Gamma
$--module with action of $\Gamma $, then the invariant $M^\Gamma $
is a direct summand of $M$.  Indeed, the map
\begin{equation}
m\longmapsto \frac{1}{|\Gamma |}\sum _{\gamma \in \Gamma
}m^\gamma,
\end{equation}
defines a retraction $M\rightarrow M^\Gamma $.  From this (i) and
(ii) follow.  Statement (iii) follows from this discussion
combined with (\cite{H-R}, remark 2.3).
\end{proof}

It follows from the lemma and standard base change theorems for
cohomology on projective schemes, that after localizing on $\Sp
(A)$, we may assume that the groups $R^jg_*(\Omega ^i\otimes
L^{-1})$ are locally free on $\Sp (A)$ of constant rank, and
compatible with base change.  Let $T$ be the scheme--theoretic
closure of a closed point of $\Sp (A\otimes \mathbb{Q})$. The
scheme $T$ is quasi--finite and flat over $\Sp (\Z )$.  Choose a
closed point $t\in T$ at which $T/\Z $ is \'etale and for which
$\text{char}(k(t))>d$.  Then $\widetilde {\scr X}\otimes k(t)$
admits a lifting to $W_2(k(t))$, and so it suffices to prove the
theorem for the pair $(\widetilde {\scr X}\otimes k(t), \widetilde
L\otimes k(t))$.  This concludes the reduction to the positive
characteristic case.

\begin{prop} Let $k$ be a field, $\scr X/k$ a tame smooth Deligne--Mumford
stack, $\pi :\scr X\rightarrow X$ its coarse moduli space, and
assume that  $\pi $ is flat. If $L$ is an invertible sheaf on
$\scr X$ such that some tensor power descends to an ample sheaf on
$X$, then there exists an integer $N_0$ such that $H^j(\scr X,
\Omega ^i_{\scr X/k}\otimes L^{-N}) = 0$ for $j<d$, all $i$, and
all $N\geq N_0$.
\end{prop}
\begin{proof}
Let $b$ be an integer such that $L^{\otimes b} = \pi ^*M$ for some
ample line bundle $M$ on $X$.

Because  $X$ is Cohen--Macaulay, $X$ has a dualizing sheaf $\omega
_X$. Since $M$ is ample, there exists an integer $l_0$ such that
for all $i$, $j<d$, $0\leq r< b$, and $l\geq l_0$
\begin{equation}
H^{d-j}(X, \omega _X\otimes (\pi _*(\Omega ^i_{\scr X}\otimes
L^{-r}))^*\otimes M^l) = 0,
\end{equation}
where $(\pi _*(\Omega ^i_{\scr X}\otimes L^{-r}))^*$ denotes the
dual of the locally free (by (\ref{2:lem} (ii))) sheaf $\pi
_*(\Omega ^i_{\scr X}\otimes L^{-r}).$  We claim that $N_0 = bl_0$
works in the proposition.  Indeed, for any $N\geq N_0$, write $N =
r+bl$, where $0\leq r< b$ and $l\geq l_0$.  Then applying
(\ref{2:lem} (i)) and Serre duality, we have
\begin{equation}
H^j(\scr X, \Omega ^i_{\scr X/k}\otimes L^{-N}) = H^{d-j}(X,
\omega _X\otimes (\pi _*(\Omega ^i_{\scr X}\otimes
L^{-r}))^*\otimes M^l)^* = 0
\end{equation}
for $j<n$.
\end{proof}

\section{Locally free log structures}

The stack $\scr X$ in (\ref{thm3}) will be constructed as a moduli
space classifying certain log structures (in the sense of Fontaine
and Illusie (\cite{Kato})).  Thus before giving the construction
and proving (\ref{thm3}) in the next section, we gather together a
few basic facts we need about log structures.  For the basic
definitions and results concerning log structures, we refer to the
first 2 sections of (\cite{Kato}).

Recall that if $P$ is a unit free monoid, then an element $p\in P$
is called irreducible if for every $p_1, p_2\in P$ for which $p =
p_1 + p_2$, either $p_1 = p$ or $p_2 = p$.  We denote the set of
irreducible elements in $P$ by $\text{Irr}(P)$.

\begin{defn}\label{monoiddefs} (i). A monoid $F$ is called \emph{free} if it is
isomorphic to $\mathbb{N}^r$ for some $r$.  The integer $r$ is
uniquely determined and is called the \emph{rank} of $F$.

\noindent (ii). A morphism $\varphi :F_1\rightarrow F_2$ between
free monoids is called \emph{simple} if $F_1$ and $F_2$ have the
same rank,  $\varphi $ is injective, and if for every irreducible
element $f_1\in F_1$ there exists an irreducible element $f_2\in
F_2$ and an integer $n$ such that $nf_2 = \varphi (f_1)$.

\noindent (iii).  A \emph{locally free log structure} on a scheme
$X$ is a fine log structure $\mc M$ on $X$ such that for every
geometric $\bar x\rightarrow X$ the monoid $\overline {\mc
M}_{\bar x}$ is free.

\noindent (iv). A morphism $\varphi :\mc M_1\rightarrow \mc M_2$
between locally free log structures on a scheme $X$ is called
\emph{simple} if for every geometric point $\bar x\rightarrow X$,
the map $\varphi _{\bar x}:\overline {\mc M}_{1, \bar
x}\rightarrow \overline {\mc M}_{2, \bar x}$ is simple in the
sense of (ii) above.
\end{defn}

\begin{rem}\label{3:rem} If $\varphi :F_1\rightarrow F_2$ is a simple morphism
as in (\ref{monoiddefs} (ii)), then for each $f_1\in
\text{Irr}(F_1)$, the element $f_2\in \text{Irr}(F_2)$ for which
there exists an integer $n$ such that $nf_2 = \varphi (f_1)$ is
uniquely determined.  From this it follows that $\varphi $ induces
a canonical bijection $\text{Irr}(F_1)\rightarrow
\text{Irr}(F_2).$
\end{rem}

\begin{lem}\label{3:lem2} Let $\varphi :\mc M_1\rightarrow \mc M_2$ be a simple
morphism of log structures on a scheme $X$, and let $\beta _1
:\mathbb{N}^r\rightarrow \mc M_1$ be a chart for $\mc M_1$. Let
$\bar x\rightarrow X$ be a geometric point, and suppose given a
surjective map $\bar \beta _2:\mathbb{N}^r\rightarrow \overline
{\mc M}_{2, \bar x}$ and integers $\{b_i\}_{i=1}^r$ prime to the
characteristic of $X$  such that the diagram
\begin{equation}\label{commdiag1}
\begin{CD}
\mathbb{N}^r@>\rho >> \mathbb{N}^r\\
@V\bar \beta _1VV @VV \bar \beta _2V\\
\overline {\mc M}_{1, \bar x}@>\bar \varphi _{\bar x}>> \overline
{\mc M}_{2, \bar x}
\end{CD}
\end{equation}
commutes, where $\rho := \oplus (\times b_i)$.  Then in some
\'etale neighborhood of $\bar x$, there exists a chart $\beta
_2:\mathbb{N}^r\rightarrow \mc M_2$ lifting $\bar \beta _2$ such
that the diagram
\begin{equation}
\begin{CD}
\mathbb{N}^r@>\rho >> \mathbb{N}^r\\
@V\beta _1VV @VV \beta _2V\\
{\mc M}_1@>\varphi >> {\mc M}_{2}
\end{CD}
\end{equation}
commutes.
\end{lem}
\begin{proof} By (\cite{Kato}, 2.10), it suffices to find a lifting $\beta
_2:\mathbb{N}^r\rightarrow \mc M_{2, \bar x}$ lifting $\bar \beta
_2$ such that the diagram
\begin{equation}\label{diag3}
\begin{CD}
\mathbb{N}^r@>\rho >> \mathbb{N}^r\\
@V\beta _1VV @VV \beta _2V\\
{\mc M}_{1, \bar x}@>\varphi >> {\mc M}_{2, \bar x}
\end{CD}
\end{equation}
commutes.  For this, choose for each $1\leq i\leq r$ a lifting
$m_i\in \mc M_{2, \bar x}$ of $\bar \beta _2(e_i)$, where $e_i$
denotes the $i$--th standard generator for $\mathbb{N}^r$.  By the
commutativity of (\ref{commdiag1}), there exists a unit $u_i\in
\mc O_{X, \bar x}^*$ such that $\lambda (u_i) + b_im_i= \varphi
(\beta _1 (e_i))$, where $\lambda :\mc O_{X, x}^*\rightarrow \mc
M_{X, \bar x}$ denotes the canonical inclusion.  Since by
assumption $b_i$ is invertible on $X$, there exists a unit $v_i\in
\mc O_{X, \bar x}^*$ such that $v_i^{b_i} = u_i$.  Replacing $m_i$
by $m_i+\lambda (v_i)$, we can therefore find units $m_i\in \mc
M_{2, \bar x}$ such that $b_im_i = \varphi (\beta _1(e_i))$.  That
is, a map $\beta _2:\mathbb{N}^r\rightarrow \mc M_{2, \bar x}$
such that (\ref{diag3}) commutes.
\end{proof}

\section{Construction of $\scr X$}

\begin{thm}\label{thm3}
Let $X/k$ be a smooth scheme over a field $k$,  $D = \bigcup
_{i\in I}D_i\subset X$ a divisor with normal crossings, and
$\{b_i\}_{i\in I}$ a collection of positive integers prime to
$\text{\rm char}(k)$. Attached to this data is a canonical pair
$(\scr X, \widetilde D = \bigcup _{i\in I}\widetilde D_i)$,
consisting of a smooth tame Deligne--Mumford stack $\pi :\scr
X\rightarrow X$ together with a normal crossings divisor
$\widetilde D\subset \scr X$ such that:

\noindent {\rm (i). } The map $\pi $ is finite and flat,
identifies $X$ with the coarse moduli space of $\scr X$, and is an
isomorphism over $X-D$.

\noindent {\rm (ii). } The pullback $\pi ^*\mc O_X(-D_i) $ is
equal to $\mc O_{\scr X}(-b_i\widetilde D_i)$ as a subsheaf of
$\mc O_{\scr X}$.

\noindent {\rm (iii).} For any collection $\{a_i\}_{i\in I}$ of
positive integers with $b_i\nmid a_i$, there is a natural
quasi--isomorphism
\begin{equation}
\Omega _{X/k}^i(\text{\rm log } D)\otimes \mc O_X(-\sum _i \ru
{a_i/b_i}D_i)\longrightarrow R\pi _*(\Omega ^i_{\scr X/k}\otimes
\mc O_{\scr X}(-\sum _ia_i\widetilde D_i)),
\end{equation}
where for $q\in \mathbb{Q}$, we denote by $\ru {q}:= \text{\rm
inf}\{p\in \Z|p\geq q\}.$

\noindent {\rm (iv).} Let $\mc M_X$ denote the log structure on
$X$ associated to the divisor $D$, and suppose $A\rightarrow k$ is
a surjective map from an Artin local ring $A$.  Then associated to
any  log smooth lifting $(\widetilde X, \mc M_{\widetilde X})$ of
$(X, \mc M_X)$ to $\Sp (A)$, there is a canonical smooth lifting
$\widetilde {\scr X}$ of $\scr X$ to $\Sp (A)$.

\noindent {\rm (v).} If $f:Y\rightarrow X$ is a finite flat
morphism from a smooth scheme $Y/k$ for which $f^*D$  is integral
with support a divisor with normal crossings, then there exists a
unique map $Y\rightarrow \scr X$ over $X$.
\end{thm}
\begin{proof} Let $\mc M_X$ be the log structure associated to the divisor
$D$ (\cite{Kato}, 1.5 (i)). The log structure  $\mc M_X$ is
locally free by (\cite{Kato}, 2.5 (i)). If $\bar x\rightarrow X$
is a geometric point, let $C(\bar x)$ denote the irreducible
components of the inverse image of $D$ in $\text{Spec}(\mc O_{X,
\bar x})$.

\begin{lem}\label{4:lem} For every geometric point $\bar x\rightarrow X$ there
is a canonical isomorphism $\overline {\mc M}_{X, \bar x}\simeq
\bigoplus _{C(\bar x)}\mathbb{N}$.
\end{lem}
\begin{proof} It suffices to establish a bijection
$\text{Irr}(\overline {\mc M}_{X, \bar x})\rightarrow C(\bar x)$.
We do so by sending an irreducible element $\bar m$ to the
component defined by the ideal $(\alpha (m))$, where $m\in \mc
M_{X, \bar x}$ is any lifting of $m$ and $\alpha :\mc M_{X, \bar
x}\rightarrow \mc O_{X, \bar x}$ is the logarithm map.  That this
is well--defined and a bijection can be verified after replacing
$X$ by an \'etale cover.  Hence it suffices to consider the case
when $x$ is the point $x_1 = \cdots = x_n = 0$ on $X = \Sp (k[x_1,
\dots , x_n])$ and $D = Z(x_1\cdots x_r),$ for some $n$ and $r$.
In this case the result is clear.
\end{proof}

If $f:Y\rightarrow X$ is any morphism of schemes, and $f^*\mc
M_X\rightarrow \mc M_Y$ is a simple morphism on $Y$, then for
every geometric point $\bar y\rightarrow Y$ with image $\bar x =
f(\bar y)$, the map
\begin{equation}
\overline {\mc M}_{X, \bar x}\longrightarrow \overline {\mc M}_{Y,
\bar y}
\end{equation}
has by (\ref{4:lem}) and (\ref{3:rem}) a canonical decomposition
\begin{equation}
\oplus (\times c_i):\bigoplus _{C(\bar
x)}\mathbb{N}\longrightarrow \bigoplus _{C(\bar x)}\mathbb{N},
\end{equation}
where $\{c_i\}$ is a collection of positive integers indexed by
$C(\bar x)$.

Define $\pi :\scr X\rightarrow X$ to be the fibered category whose
fiber $\scr X(Y)$ over $f:Y\rightarrow X$ is the groupoid of
simple morphisms of log structures $\varphi :f^*\mc M_X\rightarrow
\mc M$, such that for each geometric point $\bar y\rightarrow Y$
with $\bar x = f(\bar y)$, the integer $c_i$ associated to a
component $C_i\in C(\bar x)$ is equal to the integer $b_i$
attached (by our assumptions) to the component $D_i\subset X$
containing the image of $C_i$ under the natural map $\Sp (\mc
O_{X, \bar x})\rightarrow X$.  A morphism
\begin{equation}
\rho :(\varphi _1:f^*\mc M_X\rightarrow \mc M_1)\longrightarrow
(\varphi _2:f^*\mc M_X\rightarrow \mc M_2)
\end{equation}
in $\scr X(Y)$ is an isomorphism of log structures $\rho :\mc
M_1\rightarrow \mc M_2$ such that $\rho \circ \varphi _1 = \varphi
_2$.  With the natural notion of pullback, $\scr X$ is a fibered
category.  Moreover, because $\scr X$ classifies \'etale sheaves,
$\scr X$ is a stack with respect to the \'etale topology.

\begin{lem}\label{4:prop} The stack $\scr X$ is a Deligne--Mumford stack.  If
$X = \Sp (k[x_1, \dots , x_n])$ and $D_i = Z(x_i)$ for $1\leq
i\leq r$, then $\scr X$ is canonically isomorphic to
\begin{equation}\label{stack2}
[\Sp (k[y_1, \dots , y_n])/\mu _{b_1}\times \cdots \times \mu
_{b_r}],
\end{equation}
where $k[y_1, \dots , y_n]$ is viewed as a $k[x_1, \dots ,
x_n]$--algebra via the map induced by
\begin{equation}
x_i\mapsto \left \{\begin{array}{cl}
    y_i^{b_i} & \mbox{if $i\leq r$}\\
    y_i& \mbox{if $i>r,$}
    \end{array}\right.
\end{equation}
and $\mu _{b_1}\times \cdots \times \mu _{b_r}$ acts by
\begin{equation}
(u_1, \dots , u_r)\cdot y_i=  \left \{\begin{array}{cl}
    u_iy_i& \mbox{if $i\leq r$}\\
    y_i& \mbox{if $i>r.$}
    \end{array}\right.
\end{equation}
\end{lem}
\begin{proof} Since $\scr X$ is a stack with respect to the
\'etale topology, the assertion that $\scr X$ is a
Deligne--Mumford stack can be verified \'etale locally on $X$.
Thus it suffices to prove the second statement.  Let us
temporarily denote the stack (\ref{stack2}) by $\scr X'$.  By
(\cite{Ols}, 5.20), $\scr X'$ is canonically isomorphic to the
stack which to any $f:Y\rightarrow X$ associates the groupoid of
morphism of log structures $f^*\mc M_X\rightarrow \mc M$ together
with a map $\bar \beta ':\mathbb{N}^r\rightarrow \overline {\mc
M}$ which \'etale locally lifts to a chart, such that the diagram
diagram
\begin{equation}\label{diag4}
\begin{CD}
\mathbb{N}^r@>\oplus (\times b_i)>> \mathbb{N}^r\\
@V\bar \beta VV @VV\bar \beta 'V\\
f^{-1}\overline {\mc M}_X@>>> \overline {\mc M}
\end{CD}
\end{equation}
commutes, where $\bar \beta $ denotes the map induced by the
standard chart on $\mc M_X$.  It follows that there is a natural
functor $F:\scr X'\rightarrow \scr X$ which simply forgets the
additional data of the map $\bar \beta '$.  By (\ref{3:lem2}),
every object of $\scr X(Y)$ is \'etale locally in the essential
image of $F$, so to prove that $F$ is an equivalence, it suffices
to prove that it is fully faithful.  This amounts to the following
statement. Given two objects $(\varphi _i:f^*\mc M_X\rightarrow
\mc M_i, \beta _i')$ of $\scr X'(Y)$, any isomorphism $\rho :\mc
M_1\rightarrow \mc M_2$, such that $\rho \circ \varphi _1 =
\varphi _2$, has the property that $\bar \rho \circ \bar \beta _1'
= \bar \beta _2'$, where $\bar \rho :\overline {\mc
M}_1\rightarrow \overline {\mc M}_2$ denotes the map induced by
$\rho $. But this is clear, for since $\overline {\mc M}_2$ has no
units, two maps $\mathbb{N}^r\rightarrow \overline {\mc M}_2$ are
equal if and only if their restrictions to the image of $\oplus
(\times b_i)$ are equal.  Thus $F$ is fully faithful, and the
lemma follows.
\end{proof}

We can now complete the proof of (\ref{thm3}).  We define
$\widetilde D_i$ to be the unique reduced closed substack of $\scr
X$ whose support is equal to the support of the closed substack
defined by $\pi ^*\mc O_X(-D_i)$.  By the local description of
$\scr X$ given above, $\widetilde D := \bigcup \widetilde D_i$ is
a reduced divisor with normal crossings, and (\ref{thm3} (ii))
holds.  Also the local description of $\scr X$ in (\ref{4:prop})
implies that $X$ is the coarse moduli space of $\scr X$ (see for
example (\cite{A-V}, 2.3.3)).

To see (\ref{thm3} (iii)), note that by (\ref{2:lem} (i)), $R^i\pi
_*(\Omega ^j_{\scr X/k}\otimes  \mc O_{\scr X}(-\sum
_ia_i\widetilde D_i)) = 0$ for $i>0$. Moreover, $\Omega ^j_{\scr
X/k}(\text{log } \widetilde D) = \pi ^*\Omega ^j_{X/k}(\text{log }
D)$, so by the projection formula we have
\begin{equation}
R\pi _*(\Omega ^i_{\scr X/k}\otimes  \mc O_{\scr X}(-\sum
_ia_i\widetilde D_i))\simeq \Omega ^i_{X/k}(\text{log } D)\otimes
\pi _*\mc O_{\scr X}(-\sum _i(a_i+1)\widetilde D_i)).
\end{equation}
Since $b_i\nmid a_i$, $\ru {(a_i+1)/b_i} = \ru {a_i/b_i}$, so it
suffices to exhibit an equality of sheaves of ideals
\begin{equation}
\pi _*\mc O_{\scr X}(-\sum _i(a_i)\widetilde D_i) = \mc O_X(-\sum
_i \ru {a_i/b_i}D_i),
\end{equation}
for any collection of integers $\{a_i\}$. For this we may again
work \'etale locally on $X$, and so need only consider $X$ and
$\scr X$ as in (\ref{4:prop}). In this case, $\pi _*\mc O_{\scr
X}(-\sum _ia_i\widetilde D_i)$ is simply the $k[x_1, \dots ,
x_n]$--submodule of $(y_1^{a_1}\cdots y_r^{a_r})$ invariant under
the action of $\mu _{b_1}\times \cdots \times \mu _{b_r}$. If
$l_i$ denotes the smallest multiple of $b_i$ such that $a_i\leq
l_i$, then this submodule of invariants is $(y_1^{l_1}\cdots
y_r^{l_r})$. Since $l_i = \ru {a_i/b_i} \ b_i$, this implies
(\ref{thm3} (iii)).

As for (\ref{thm3} (iv)), let $\widetilde {\scr X}$ be the fibered
category on the category of $\widetilde X$--schemes, which to any
$f:Y\rightarrow \widetilde X$ associates the groupoid of morphisms
of log structures $f^*\mc M_{\widetilde X}\rightarrow \mc M_Y$
which when pulled back to $Y\otimes _Ak$ defines an object of
$\scr X(Y\otimes _Ak)$.  Since $(\widetilde X, \mc M_{\widetilde
X})$ is log smooth over $\Sp (A)$, it follows from the structure
theorem for log smooth morphisms (\cite{Kato}, 3.4) that \'etale
locally on $\widetilde X$, the log scheme $(\widetilde X, \mc
M_{\widetilde X})$ is isomorphic to $\Sp (A[x_1, \dots , x_n])$
with log structure defined by
\begin{equation}
\mathbb{N}^r\longrightarrow A[x_1, \dots , x_n], \ \ \ (i_1, \dots
, i_r)\mapsto x_1^{i_1}\cdots x_r^{i_r},
\end{equation}
for some $n$ and $r$.  From this and the same argument used in the
proof of (\ref{4:prop}), it follows that $\widetilde {\scr X}$ is
a smooth lifting of $\scr X$.

To prove (\ref{thm3} (v)), we have to attach to any
$f:Y\rightarrow X$ as in the theorem, a unique (up to unique
isomorphism) object $f^*\mc M_X\rightarrow \mc M$ of $\scr X(Y)$.
Since $f$ is flat, $f^*\mc M_X$ is a subsheaf of $\mc O_Y$.  From
this it follows that $\mc M$ is also a subsheaf of $\mc O_Y$.
Moreover, we can completely characterize the sheaf $\mc M$ as
follows.  Etale locally on $X$, we can find a chart
$\mathbb{N}^r\rightarrow \mc  M_X$ for the log structure on $X$.
For each generator $e_i\in \mathbb{N}^r$, let $f_i$ be its image
in $\mc O_X$ and  let $b_i$ denote the corresponding integer
attached to the component of $D$ defined by $f_i$.  By assumption,
there exists \'etale locally on $Y$ an element $g_i\in \mc O_Y,$
well defined up to multiplication by a root of unity,  with
$g_i^{b_i} = f_i$. The log structure $\mc M$ is then the subsheaf
of monoids of $\mc O_Y$ generated by the $g_i$ and $\mc O_Y^*$.
This also proves the existence.
\end{proof}

\section{The Kawamata--Viehweg vanishing theorem}

 Let $k$ be a  field,  $X/k$ be a smooth projective variety of pure dimension
 $d$, and
 $E = \sum
_ia_iD_i\in \text{Div}(X)\otimes \mathbb{Q}$ an ample divisor.
Define the \emph{integral part} $[E]$, the \emph{round up} $\ru
E$, and the \emph{fractional part} $\langle E\rangle$  by
\begin{align}
[E] & := \sum [a_i]D_i\\
\ru E & := \sum \ru {a_i} D_i = -[-E]\\ \langle E \rangle & :=
\sum \langle a_i\rangle D_i = E-[-E],
\end{align}
where for $r\in \mathbb{R}$ we write $[r] := \text{max}\{t\in
\mathbb{Z}|t\leq r\}.$  We suppose that the support of $\fp E$ is
a divisor $D$ with normal crossings, and denote by $\mc M_X$ the
associated log structure on $X$ (\cite{Kato}, 1.5 (i)).  The
following version of the Kawamata--Viehweg vanishing theorem can
be found in (\cite{K-M-M}, 1--2--2):
\begin{thm}\label{thm1} {\rm (i).} Suppose $\text{\rm char}(k) = 0$.
Denote by $\Omega _X^1(\text{\rm log } D)$ the sheaf of
differentials on $X$ with log poles along $D$, and by $\Omega
_X^i(\text{\rm log } D)$ its $i$--th exterior power. Then
\begin{equation}
H^j(X, \Omega _X^i(\text{\rm log } D)\otimes \mc O_X(-\ru E)) = 0
\text{ for  } i+j<d.
\end{equation}

\noindent {\rm (ii).} Suppose $k$ is perfect and let $W_2(k)$
denote the reduction modulo $p^2$ of the ring of Witt vectors of
$k$. Suppose further that  the smooth log scheme $(X, \mc M_X)$
admits a log smooth lifting to $W_2(k).$ Then
\begin{equation}
H^j(X, \Omega _X^i(\text{\rm log } D)\otimes \mc O_X(-\ru E)) = 0
\text{ for  } i+j<\text{\rm inf}(d, p).
\end{equation}
\end{thm}
\begin{proof}
Write $\langle E\rangle  = \sum _i(a_i/b_i)D_i$ with $(a_i, b_i) =
1$. Note that by the openness of the ample cone, we can replace
the rational numbers $a_i/b_i$ by any other rational numbers
sufficiently close to the $a_i/b_i$ without changing the statement
of the theorem. Thus in the positive characteristic case, we may
assume that the integers $b_i$ are prime to $p$.  Let $\pi :\scr
X\rightarrow X$ be the stack associated to the data $(X, D,
\{b_i\})$ as in (\ref{thm3}). Then by (\ref{thm3} (iii)), and the
projection formula, we have
\begin{equation}
H^j(X, \Omega ^i_X(\text{log }D)\otimes \mc O_X(-\ru E)) =
H^j(\scr X, \Omega ^i_{\scr X}\otimes \mc O_{\scr X} (-\pi
^*[E]-\sum _i a_i\widetilde D_i)).
\end{equation}
Since $(\prod _ib_i)\mc O_{\scr X}(\pi ^*[E]+\sum _ia_i\widetilde
D_i)$ descends to an ample sheaf on $X$, the theorem therefore
follows from (\ref{thm2}).
\end{proof}

\begin{rem}
The above theorem plays a crucial and indispensible role in
carrying out the proofs for many key ingredients of the so-called
Minimal Model Program (MMP for short).  The call for
Kawamata--Viehweg vanishing, going beyond the classical Kodaira
vanishing, is more apparent when we have to deal with the
singularities that MMP inevitably brings in higher dimensions
($\text{dim} \geq 3$). This is usually perceived as a technical
calamity rather than an indication of any essential point.  We
wonder and/or speculate, however, that, once we have an
interpretation of Kawamata--Viehweg vanishing as a version of
Kodaira vanishing for stacks as in this note, there may be some
smooth stacks floating around behind the whole game of MMP, and
the singularities of MMP we only observe as we look at the coarse
moduli of these smooth stacks.
\end{rem}

\end{document}